\documentclass[11pt]{article}
\usepackage{graphicx}
\usepackage{amssymb}
\usepackage{epstopdf}
\def\b{\beta}
\def\C{\Gamma}
\def\c{\gamma}

\def\t{\theta}

\def\eproof{$\Box$ \medskip}

\newcommand{\CC}{\mathbb C}

\newcommand{\R}{\mathbb R}
\newcommand{\N}{\mathbb N}
\newcommand{\Hp}{{\mathbb H}^2}
\newcommand{\Hs}{{\mathbb H}^3}
\newtheorem{theorem}{Theorem}
\newtheorem{lemma}{Lemma}

\title{Higher derivatives of length functions along earthquake deformations}
\author{Martin Bridgeman\thanks{This work was partially supported by grant \#266344 from the Simons Foundation}}

\begin{document}
\maketitle

\section{Introduction}
Let $S$ be a closed surface of genus $g\geq 2$ and  $T(S)$  the associated  Teichm\"uller space of hyperbolic structures on $S$.  Given $\gamma \in \pi_1(S)$, let  $L_\gamma:T(S) \rightarrow \R$ be the associated length  function and $T_\gamma:T(S) \rightarrow \R$ the associated trace function. The functions $L_\gamma, T_\gamma$ have a simple relation given by
\begin{equation}
T_\gamma = 2\cosh(L_\gamma/2).
\label{relation}
\end{equation}

Let $\b$ be the homotopy class of a simple multicurve (i.e. a union of disjoint simple non-trivial closed curves in $S$) and $t_\b$ the vector field on $T(S)$ associated with left twist along the geodesic representative of  $\b$ (see \cite{Ker83}). In this paper, we describe a formula to calculate the  higher order derivatives of the  functions $L_\gamma, T_\gamma$ along $t_\b$. In particular we will find a formula for 
$$t_\beta^k L_\gamma = t_\beta t_\beta\ldots t_\beta L_\gamma.$$

The formulae we derive generalize formulae for the first  two derivatives due to Kerchoff (1st derivative, see \cite{Ker83}) and Wolpert (1st and 2nd derivatives, see \cite{Wol81,Wol86}).

Kerckhoff and Wolpert  both showed that the first derivative is given by
\begin{equation}
t_\beta L_\gamma = \sum_{p\in \beta' \cap \gamma'} \cos{\theta_p}
\label{1d}
\end{equation}
where $\beta', \gamma'$ are the geodesic representatives of $\beta, \gamma$ respectively and $\theta_p$ is the angle of intersection at $p \in \beta'\cap \gamma'$. Kerckhoff further generalized the formula for the case when $\beta,\gamma$ are measured laminations (see \cite{Ker83}).

In \cite{Wol86}, Wolpert  derived this formula for second derivative   
$$t_\alpha t_\beta L_\gamma = \sum_{(p,q) \in \beta'\cap\gamma' \times \alpha'\cap \gamma'} \frac{e^{l_{pq}}+e^{l_{qp}}}{2(e^{L_\gamma}-1)}\sin\theta_p\sin\theta_q +\sum_{(r,s) \in \beta' \cap \gamma' \times \beta'\cap\alpha'} \frac{e^{m_{rs}}+e^{m_{sr}}}
{2(e^{L_\beta}-1)}\sin\theta_r\sin\theta_s .$$
where  $l_{xy}$ is the length along $\gamma$ between $x,y$ (similarly $m_{xy}$ is the length along $\beta$.)

It follows from Wolpert's formula  that 
\begin{equation}
t_\b^2 L_\gamma = t_{\b}t_{\b}L_{\c} =  \sum_{p,q \in \b' \cap \c'} \frac{e^{l_{pq}} + e^{l_{qp}}}{2(e^{L_{\c}}-1)} \sin \t_{p}\sin\t_{q}.
\label{2d}
\end{equation}

Our  formula  generalizes equations \ref{1d}, and \ref{2d} to higher derivatives. Our approach is to derive a formula for the higher derivatives of $T_\gamma$ and then use the functional relation in equation \ref{relation} to derive the formula for $L_\gamma$.

\section{Higher Derivative Formula}
We take the geodesic representatives of $\beta$ and $\gamma$. We let the geometric intersection number satisfy $i(\beta,\gamma) = n$ and we order the points of intersection $x_1,\ldots,x_n$ by choosing a base point on $\gamma$. We let $\theta_i$ be the angle of intersection of $\beta,\gamma$ at $x_i$ and $l_i$ be the length along $\gamma$ from $x_1$ to $x_i$. This gives us $n$-tuples $(l_1,\ldots, l_n)$ and $(\theta_1,\ldots,\theta_n)$.

 In order to describe the formula for the higher derivatives, we first introduce some more notation. 
 
 Given $r$, we let $P(r)$ be the set of subsets of the set $\{1,\ldots,r\}$. Then for $I \in P(r)$ will be denoted by $I = (i_1,\ldots,i_k)$ where $1 \leq i_1 < i_2 <\ldots<i_k \leq r$. We then define $\hat{I}$ to be the complementary subset. We also let $|I|$ be the cardinality of $I$. 
 
We define the alternating length  $L_I$ for $I = (i_i,\ldots, i_k)$ by
$$L_I= \sum_{j=1}^{k} (-1)^{j}l_{i_j} = -l_{i_1}+l_{i_2}-l_{i_3}-\ldots + (-1)^{k}l_{i_k}$$
We further define a signature for $I \in P(r)$.  For $I = (i_1,\ldots, i_k)$  we can consider the integers in $\{1,\ldots,r\}$ in the ordered blocks $[1,i_1],[i_1,i_2],\ldots [i_k,r]$.  We take the sum of the cardinality of the even ordered blocks. Then 
$$s(I) = (i_2-i_1 +1) + (i_4-i_3+1) +\ldots (i_k-i_{k-1} + 1)\qquad k \mbox{ even}$$
$$s(I) = (i_2-i_1 +1) + (i_4-i_3 +1) +\ldots (r-i_k +1)\qquad k \mbox{ odd}$$
For $(\theta_1,\ldots,\theta_n)$ we also define
$$\cos(\theta_I) = \prod_{j=1}^k \cos(\theta_{i_j}) = \cos(\theta_{i_1})\cos(\theta_{i_2})\ldots\cos(\theta_{i_k})$$
and similarly define $f(\theta_I)$ for $f$ a   trigonometric functions.

We let $u_j = l_j +i \theta_j$. The function $F_r$ is given by 
$$F_r(u_1,\ldots,u_r, L) = \sum_{I \in P(r), |I| \mbox{ even}} (-1)^{s(I)}\sin(\theta_I)\cos(\theta_{\hat{I}})\left(e^{L/2-L_I}+ (-1)^{r}e^{L_I-L/2}\right).$$
or equivalently 
$$F_r(u_1,\ldots,u_r, L) = \sum_{I \in P(r), |I| \mbox{ even}} (-1)^{s(I)}2\sin(\theta_I)\cos(\theta_{\hat{I}})\cosh(L/2-L_I)  $$
for  $r$ even and
$$F_r(u_1,\ldots,u_r,L) = \sum_{I \in P(r), |I| \mbox{ even}} (-1)^{s(I)}2\sin(\theta_I)\cos(\theta_{\hat{I}})\sinh(L/2-L_I)$$
for $r$ odd.

We let $C(n,r)$ be the set of  the subsets of size $r$ of the set $\{1,2,\ldots,n\}$. It is given by
$$C(n,r) =  \left\{ I =(i_1,i_2,\ldots,i_r)\ \left|\  1\leq i_1 < i_2 < \ldots < i_r \leq n \right\}\right.$$

Given $m \in \N$, we let $[m]$ be the parity of $m$, i.e. $[m] = 0$ if $m$ is even, and $[m] = 1$ if $m$ is odd.

\begin{theorem}
Let $\beta$ be a homotopy class of a simple closed multicurve and $\gamma$ a homotopy class of non-trivial closed curve. Let the geometric intersection number $i(\beta,\gamma) = n$.  Then 
$$t_\beta^k T_\gamma = \frac{1}{2^k}\sum_{\stackrel{r = 0}{[r]=[k]}}^{k} B_{n,k,r}  \sum_{I \in C(n,r)}  F_{r}(u_{i_1},\ldots, u_{i_r}, L_\gamma) $$ 
where $B_{n,k,r}$ are constants described below.
\end{theorem}

The first two equations correspond to  formulae \ref{1d} and \ref{2d} for the derivatives of length. We use the above, to derive the next case as an example.

{\bf Third Derivative:}
We use the above formula to calculate the formula for the third derivative.
$$t_\beta^{3}T_\gamma= \frac{1}{8} ((6n-4)\sinh(L_\gamma/2)\sum_{i=1}^n \cos(\theta_i) +12(\sum_{i<j <k}\sinh(L_\gamma/2)\cos(\t_i) \cos(\t_j) \cos(\t_k) $$
$$  + \sinh(L_\gamma/2-l_{ij})\sin(\t_i) \sin(\t_j)\cos(\t_k)  - \sinh(L_\gamma/2-l_{ik})  \sin(\t_i)\cos(\t_j)\sin(\t_k) $$
$$\left.+\sinh(L_\gamma/2-l_{jk})  \cos(\t_i)\sin(\t_j) \sin(\t_k)\right)$$

\subsection{Constants $B_{n,k,r}$}

We let $P(k,n)$ be the collection of partitions of $k$ into $n$ ordered nonnegative integers, i.e.
$$P(k,n) = \left\{ p =(p_1,p_2,\ldots,p_n) \in \N_0^{\ n} \ \left|\  \sum_{i=1}^n p_i = k \right\}\right.$$

For $p \in P(k,n)$, we define $[p] = ([p_1],\ldots,[p_n])$ where $[n]$ is the parity of $n$. We let $|p| = [p_1]+\ldots[p_n]$. Then $[p]$  is an n-tuple of 0's and 1's with exactly $|p|$ 1's. 

Given $p \in P(k,n)$ we define $B(p)$ as a sum of multnomials given by
$$B(p) = \sum_{q \in P(k,n), [q] = [p]}  \left( \begin{array}{c}
k\\
q\end{array}
\right).
$$

It is easy to see that $B(p)$ only depends on $n,k$ and $r =|p|$. We therefore define 
$$B_{n,k,r} = B(p)\qquad \mbox{ for some $p$ with $|p| = r$}$$
 In particular if we let $p_r = (1,1,\ldots,1,0,\ldots,0) \in P(k,n)$, of $r$ 1's followed by $(n-r)$ 0's, we have
$$B_{n,k,r} = \sum_{p \in P(k,n), [p] = [p_r]}  \left( \begin{array}{c}
k\\
p\end{array}
\right).
$$
A simple calculation gives 
$$B_{n,k,k} =   \left( \begin{array}{c}
k\\
p_k\end{array}
\right) =  \left( \begin{array}{c}
k\\
1,1,1\dots,0,0,\ldots 0\end{array}
\right) = k!
$$

\section{Twist Deformation}
We consider $T(S)$ as the fuchsian locus  of the associated quasifuchsian space  $QF(S)$.  Let $X \in T(S)$ and 
$X = \Hp/\C$ where  $\C$ is a subgroup of $PSL(2,{\bf \CC})$ acting on upper half space  $\Hs =\{(u,v,w) \in \R^3\ |\  w > 0\}$ fixing the hyperbolic plane $\Hp = \{(u, 0 ,w) | w > 0\}$.  Let $\C_{z}$ be the subgroup of $PSL(2,\CC)$ obtained by complex shear-bend along $\b$ by amount $z = s + it$, i.e. left shear by amount $s$ followed by bend of $t$. Then for small $z$, $X_z = \Hs/\C_z \in QF(S)$. In the terminology of Epstin-Marden this is a quake-bend deformation. See  II.3 of \cite{EM05} for details on quake bend deformations and  II.3.9 for a detailed discussion of derivatives of length along quakbend deformations.

 Let $\gamma \in \C$ be a hyperbolic element and let $\gamma(z) \in \C_{z}$ be the element of the deformed group corresponding to $\gamma$ and  $L(z)$  the complex translation length of $\gamma(z)$.  To see how $\gamma$ is deformed, by conjugating, we  assume that $\gamma$ has axis the geodesic $g$ with endpoints $0,\infty \in \hat{\CC}$ and is given by
$$\gamma = \left(\begin{array}{cc}
\lambda & 0 \\
0 & 1/\lambda\end{array}\right)
\mbox{ with } \lambda = e^{L/2} \mbox{ where $L>0$ is the translation length of $\gamma$}.$$
We consider the lifts of $\beta$ which intersect the axis $g$ of $\gamma$ and normalize to have  a lift of $\b$ labelled  $\b_{1}$ which intersects axis $g$ at height $1$. We enumerate all other lifts by the order of the height of their  intersection point with $g$ starting with the intersection point of $\b_{1}$. Let $n$ be such that $\gamma \b_{1} = \b_{n+1}$. Let $R_{i}(z)$ be the M\"obius transformation corresponding to a complex bend about $\b_{i}$ of  $z$. Then under the complex bend about $\b$,  $\gamma(z)$ given by
$$\gamma(z) = R_{1}(z)R_{2}(z)\ldots R_{n}(z)\gamma.$$
A similar description of the deformation of an element in the punctured surface case can be given in terms of shearing coordinates (see \cite{Chek07} for details).

Taking traces we have
$$T(z) = Tr(R_{1}(z)R_{2}(z)\ldots R_{n}(z)\gamma) = 2\cosh(L(z)/2).$$
We can find the derivatives of $L(z)$ by differentiating this formula repeatedly. The final formula is obtained by applying symmetry relations on the derivatives and some elementary combinatorics.

We note that both $T(z)$ and $L(z)$ are holomorphic in $z$. Differentiating in the real direction we have
$$ t_\b^k L_\gamma = \frac{d^kL}{dz^k}(0) = L^{(k)}(0) $$
Also if we let $b_\b$ be the vector field on $T(S)$ given by pure bending along $\b$, then we have by analyticity of $L(z)$
$$b_\b^k L_\gamma =  i^k L^{(k)}(0) = (it_\b)^k L_\gamma.$$
This corresponds  to the observation that $b_\b = J.t_\b$ where $J$ is the complex structure on $QF(S)$ (see \cite{Bon96}). 

\begin{figure}[htbp] 
   \centering
   \includegraphics[width=4in]{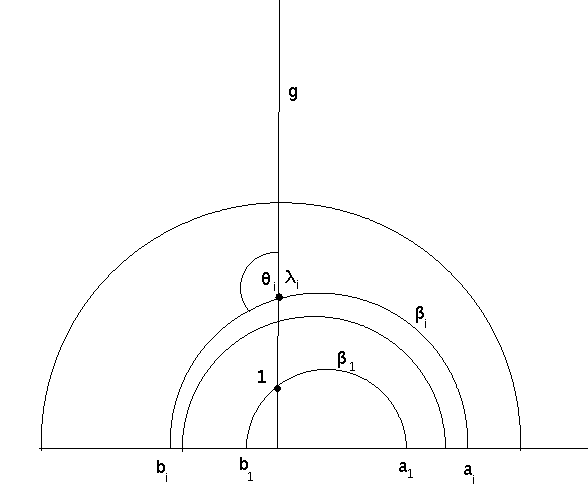} 
   \caption{Lift of $\gamma$}
   \label{twist}
\end{figure}

\subsection{Derivation of First Two Derivatives}
We now  calculate the first two derivatives and recover Wolpert's formulae. By the product rule we have
\begin{eqnarray}
T'(0) = \sum_{i = 1}^{n} Tr(R'_{i}(0)\gamma) \qquad
T''(0) =  \sum_{i=1}^{n}Tr(R''_{i}(0)\gamma) +  2\sum_{\stackrel{i,j =1}{i< j}}^{n} Tr(R'_{i}(0)R'_{j}(0)\gamma)
\label{traceeqn}
\end{eqnarray}
We now describe $R_{i}(z)$. Let $\b_{i}$ have endpoints $a_{i}, b_{i} \in {\bf R}$ where $a_{i} > 0$ and $b_{i} < 0$. We  let $\lambda_i$ be the height at which $\beta_i$ intersects $g$. We orient $\beta_i$ from $a_i$ to $b_i$ and $g$ from $0$ to $\infty$ and let $\theta_{i}$ be the angle $\b_{i}$ makes with side $g$ with respect to these orientations (see figure \ref{twist}). 

Then
$$\lambda_i = \sqrt{-a_{i}b_{i}} \qquad \cos{\theta_{i}} = -\left(\frac{a_{i}+b_{i}}{a_{i}-b_{i}}\right)\qquad 
\sin{\theta_{i}} = \frac{2\sqrt{-a_{i}b_{i}}}{a_{i}-b_{i}}.$$ 
As $\beta_1$ intersects at height 1, the distance $l_i$ between the intersection points of $\beta_1$ and $\beta_i$ is given by $e^{l_i} = \lambda_i.$ Then we let $f_{i} \in SL(2,\R)$ acting on the upper-half space by $f_i(z) = (z-a_{i})/(z-b_{i})$ and let $S(z) = f_{i} R_{i}(z) f_{i}^{-1}$. Then $S_{i}(z)$ is the complex translation given by 
$$S(z) = \left( \begin{array}{cc}
e^{z/2} & 0\\
0 & e^{-z/2}\end{array}
\right).$$
Thus as $R_{i}(z) = f_{i}^{-1}S(z) f_{i}$. Taking derivatives  we have $R'_i(0) =  f_{i}^{-1}S'(0) f_{i}$ and
$$R'_{i}(0) = \frac{1}{a_{i}-b_{i}}\left( \begin{array}{cc}
-b_{i} & a_{i}\\
-1 & 1 \end{array}\right) \left( \begin{array}{cc}
1/2 & 0\\
0 & -1/2\end{array}\right) \left( \begin{array}{cc}
1 & -a_{i}\\
1 & -b_{i}\end{array}
\right) = \frac{1}{2(a_{i}-b_{i})}\left( \begin{array}{cc}
-(a_{i}+b_{i}) & 2a_{i}b_{i}\\
-2 & a_{i} + b_{i}\end{array}
\right)$$
Therefore
$$R'_{i}(0) = \frac{1}{2}\left( \begin{array}{cc}
\cos\theta_i & -e^{l_i}\sin\theta_i\\
-e^{-l_i}\sin\theta_i & -\cos\theta_i\end{array}
\right)$$ 
Also as $S''(0) = \frac{1}{4}I$ we have $R''_i(0) =  \frac{1}{4}I$.
Using this we have that
$$Tr(R'_{i}(0)\gamma) = Tr\left(\frac{1}{2}\left( \begin{array}{cc}
\cos\theta_i & -e^{l_i}\sin\theta_i\\
-e^{-l_i}\sin\theta_i & -\cos\theta_i\end{array}
\right)\left( \begin{array}{cc}
e^{L/2} & 0\\
0 & e^{-L/2}\end{array}
\right)\right) = \sinh(L/2)\cos\theta_i  $$
$$Tr(R'_{i}(0)R'_{j}(0)\gamma) = Tr\left(\frac{1}{4}\left( \begin{array}{cc}
\cos\theta_i & -e^{l_i}\sin\theta_i\\
-e^{-l_i}\sin\theta_i & -\cos\theta_i\end{array}
\right)\left( \begin{array}{cc}
\cos\theta_j & -e^{l_j}\sin\theta_j\\
-e^{-l_j}\sin\theta_j & -\cos\theta_j\end{array}
\right)\left( \begin{array}{cc}
e^{L/2} & 0\\
0 & e^{-L/2}\end{array}
\right)\right).$$
\begin{equation} 
= \frac{1}{4}\left(\cos\t_{i} \cos\t_{j}(e^{L/2} + e^{-L/2}) + \sin\t_{i}\sin\t_{j}\left(e^{L/2+l_i-l_j} + e^{-(L/2+l_i-l_j)}\right)\right)
\end{equation}

Let $l_{ij}$ be the distance along $\gamma$ from $\beta_i$ to $\beta_j$ with respect to the  orientation of $\gamma$. Then for $i <j$ we have $l_{ij} = l_j-l_i$,  and $l_{ji} = L-l_{ij}$ for $i> j$.
\begin{equation}
Tr(R'_{i}(0)R'_{j}(0)\gamma) = \frac{1}{2}(\cos\t_{i} \cos\t_{j}\cosh(L/2) + \sin\t_{i}\sin\t_{j}\cosh(L/2-l_{ij})).
\label{11d}
\end{equation}
Combining these we obtain the first two derivatives of $T_\gamma$.
$$T'(0) =  \sinh(L/2) \sum_{i=1}^{n}\cos\t_{i}$$
$$T''(0) = \sum_{\stackrel{i,j = 1}{i < j}}^{n}\left(\cos\t_{i} \cos\t_{j} \cosh(L/2)  + \sin\t_{i}\sin\t_{j}\cosh(L/2-l_{ij})\right) + \frac{n\cosh(L/2)}{2}$$ 

As $T(z) = 2\cosh(L(z)/2)$, then $T'(0) = \sinh(L/2)L'(0)$
giving
$$L'(0) = \sum_{i=1}^{n}\cos\t_{i}$$
Also $T''(0) = \frac{1}{2}\cosh(L/2)(L'(0))^{2} +  \sinh(L/2)L''(0).$ Therefore
$$T''(0) = \frac{\cosh(L/2)}{2}\left( n  + 2.\sum_{\stackrel{i,j = 1}{i < j}}^{n}\cos\t_{i} \cos\t_{j}\right)  + \sum_{\stackrel{i,j = 1}{i < j}} \sin\t_{i}\sin\t_{j}\cosh(L/2-l_{ij}) .$$
We have
$$n  + 2.\sum_{\stackrel{i,j = 1}{i \neq j}}^{n}\cos\t_{i} \cos\t_{j} = \left( \sum_{i=1}^{n}\cos\t_{i}\right)^{2} + \sum_{i=1}^{n}\sin^{2}{\t_{i}}$$
and 
\begin{equation}
T''(0) =   \frac{\cosh(L/2)\left(  (\sum_{i=1}^{n}\cos\t_{i})^{2} +\sum_{i=1}^{n}\sin^{2}{\t_{i}}\right)}{2} +\sum_{\stackrel{i,j = 1}{i < j}} \sin\t_{i}\sin\t_{j}\cosh(L/2-l_{ij}) .
\end{equation}
Solving for $L''(0)$ we obtain
$$L''(0) = \sum_{i=1}^{n}\frac{\sin^{2}{\t_{i}}}{2\tanh(L/2)} + \sum_{\stackrel{i,j = 1}{i < j}}^n \frac{\sin\t_{i}\sin\t_{j}\cosh(L/2-l_{ij})}{\sinh(L/2)} $$
As $l_{ii} = 0$ we can write
$$ L''(0) =  \sum_{i,j = 1}^n \frac{e^{l_{ij}-L/2} + e^{L/2 - l_{ij}}}{2(e^{L/2}-e^{-L/2})} \sin \t_{i}\sin\t_{j} = \sum_{i,j= 1}^n \frac{e^{l_{ij}} + e^{l_{ji}}}{2(e^{L}-1)} \sin \t_{i}\sin\t_{j} .$$

The above give the formulae \ref{1d} and \ref{2d} as  described.

\section{Higher Derivatives}
We now derive the formula for higher derivatives. We have the formula
$$T(z) = Tr(R_{1}(z)R_{2}(z)\ldots R_{n}(z)\gamma).$$
We let $P(k,n)$ be the collection of partitions of $k$ into $n$ ordered nonnegative integers, i.e.
$$P(k,n) = \left\{ p =(p_1,p_2,\ldots,p_n) \in \N_0^{\ n} \ \left|\  \sum_{i=1}^n p_i = k \right\}\right.$$

Then by the product rule, the kth derivative of $T$ at zero is,
$$T^{(k)}(0) = \sum_{p \in P(k,n)}  \left( \begin{array}{c}
k\\
p\end{array}
\right) Tr(R_1^{(p_1)}(0)\ldots R_n^{(p_n)}(0)\gamma)
$$

We have from above that $R_{i}(z) = f_{i}^{-1}S(z) f_{i}$ where 
$$S(z) = \left( \begin{array}{cc}
e^{z/2} & 0\\
0 & e^{-z/2}\end{array}
\right).$$
As $S^{(2)}(z) = \frac{1}{4} S(z)$, we have for $m$ even
$$R_i^{(m)}(0) = \frac{1}{2^m} I$$
and for $m$ odd we have
$$R_i^{(m)}(0) = \frac{1}{2^{m-1}} R'_i(0) =  \frac{1}{2^m}\left( \begin{array}{cc}
\cos\theta_i & -e^{l_i}\sin\theta_i\\
-e^{-l_i}\sin\theta_i & -\cos\theta_i\end{array}
\right).$$ 
Let $z = x+ i y$ and define
 $$A(z) =   \left( \begin{array}{cc}
\cos{y} & -e^{x}\sin{y}\\
-e^{-x}\sin{y} & -\cos{y}\end{array}
\right).$$

We let $u_j = l_j + i \theta_j$. Then 
$$R_j^{(p)}(0) =  \left\{ \begin{array}{cc}
\frac{1}{2^p} A(u_j)& p \mbox{ odd}\\
\frac{1}{2^p} I & p \mbox{ even}\end{array}
\right.$$
Therefore
$$T^{(k)}(0) = \frac{1}{2^k}\sum_{p \in P(k,n)}  \left( \begin{array}{c}
k\\
p\end{array}
\right) Tr(A(u_1)^{[p_1]}\ldots A(u_n)^{[p_n]}\gamma)
$$
where $[m]$ is the parity of $m$. We define
$$F_r(z_1,\ldots, z_r,L) = Tr(A(z_1)\ldots A(z_r)\gamma)$$
Therefore gathering terms we have
$$T^{(k)}(0) = \frac{1}{2^k}\sum_{r = 0}^{k} B_{n,k,r}  \sum_{1 \leq i_1 < \ldots < i_r \leq n}  F_{r}(u_{i_1},\ldots, u_{i_r}, L) $$
where $B_{n,k,r}$ are the coefficients described above. We note that we only get non-zero terms for $[r] = [k]$, so we have $B_{n,k,r} = 0$ for $[k] \neq [r]$.

We define the function 
$$G_{r}(u_1,\ldots, u_n, L) = \sum_{I \in C(n,r)}  F_{r}(u_{i_1},\ldots, u_{i_r}, L).$$
Then $G_r$ is symmetric in $(u_1,\ldots,u_n)$ and we have
$$t_\beta^k T_\gamma = \frac{1}{2^k}\sum_{\stackrel{r = 0}{[r]=[k]}}^{k} B_{n,k,r} G_{r}(u_1,\ldots, u_n, L_\gamma)$$ 

\subsection{Function $F_r$}
We now calculate the formula for $F_r$. 

\begin{lemma}
The function $F_r$ is given by 
$$F_r(u_1,\ldots,u_r, L) = \sum_{I \in P(r), |I| \mbox{ even}} (-1)^{s(I)}\sin(\theta_I)\cos(\theta_{\hat{I}})\left(e^{L/2-L_I}+ (-1)^{r}e^{L_I-L/2}\right).$$
or equivalently 
$$F_r(u_1,\ldots,u_r, L) = \sum_{I \in P(r), |I| \mbox{ even}} (-1)^{s(I)}2\sin(\theta_I)\cos(\theta_{\hat{I}})\cosh(L/2-L_I)  $$
for  $r$ even and
$$F_r(u_1,\ldots,u_r,L) = \sum_{I \in P(r), |I| \mbox{ even}} (-1)^{s(I)}2\sin(\theta_I)\cos(\theta_{\hat{I}})\sinh(L/2-L_I)$$
for $r$ odd.
\end{lemma}

{\bf Proof:}
We have $$A(u) =   \left( \begin{array}{cc}
\cos\theta & -e^{l}\sin\theta\\
-e^{-l}\sin\theta & -\cos\theta\end{array}
\right).$$
Therefore $F_r(u_1,\ldots,u_r,L) = Tr(A(u_1),\ldots,A(u_r)\gamma)$ has the form 
$$F_r(u_1,\ldots, u_r, L) =  \sum_{I \in P(r)} a_I \sin(\theta_I)\cos(\theta_{\hat{I}})$$
for some coefficients $a_I$.
Expanding we have 
\begin{eqnarray*}
F_r(u_1,\ldots,u_r,L) &= (A(u_1)\ldots A(u_r)\gamma)^1_1 +(A(u_1)\ldots A(u_r)\gamma)^2_2 \\
&=  e^{L/2}(A(u_1)\ldots A(u_r))^1_1+ e^{-L/2}(A(u_1)\ldots A(u_r))^2_2.
\end{eqnarray*}
Similarly we have
 $$(A(u_1)\ldots A(u_r))^i_j = \sum_{I \in P(r)} a^i_j(I) \sin(\theta_I)\cos(\theta_{\hat{I}})$$
 and define
  $$(A(u_1)\ldots A(u_r))^i_j (I) = a^i_j(I) \sin(\theta_I)\cos(\theta_{\hat{I}}).$$
  
We prove the lemma by induction. Given  $I = (i_1,\ldots,i_k) \in P(r)$ then $I_j = (i_{1}, i_{2},\ldots,i_{j-1}) \in P(i_j)$.

The matrix $A(u)$ has $\cos$ terms on the diagonal and $\sin$ off diagonal. As $\sin(\theta_{i_k})$ is the last $\sin$ term we have in $(A(u_1),\ldots,A(u_r))^1_1(I)$ we have 
 $$\left(A(u_1)\ldots A(u_r)\right)^1_1(I) = (A(u_1)\ldots A(u_{i_{k-1}}))^1_2(I_k) (A(u_{i_k})^2_1A(u_{i_k+1})^1_1\ldots A(u_r)^1_1$$
 $$= \cos(\theta_{i_k+1})\ldots\cos(\theta_{r})\left(-e^{-l_{i_k}}\sin(\theta_{i_k})\right)(A(u_{1})\ldots A(u_{i_{k-1}}))^1_2(I_k)$$
Now iterating, as the next sin is $\sin(\theta_{i_{k-1}})$ we have
$$(A(u_{1})\ldots A(u_{i_{k-1}}))^1_2(I_k) = (A(u_{1})\ldots A(u_{i_{k-2}}))^1_1(I_{k-1}) A^1_2(u_{i_{k-1}})A^2_2(u_{i_{k-1}+1})A^2_2(u_{i_{k-1}+2})\ldots A^2_2(u_{i_{k-1}}).$$
$$=(A(u_{1})\ldots A(u_{i_{k-2}}))^1_1(I_{k-1}) (-e^{l_{i_{k-1}}}\sin(\theta_{i_{k-1}}))(-\cos(\theta_{i_{k-1}+1}))(-\cos(\theta_{i_{k-1}+2}))\ldots(-\cos(\theta_{i_k-1})).$$
Thus we have
$$\frac{(A(u_1),\ldots,A(u_r))^1_1(I)}{(A(u_{1})\ldots A(u_{i_{k-2}}))^1_1(I_{k-1})}= $$
$$(-1)^{i_k-i_{k-1}+1}e^{l_{i_{k-1}}-l_{i_k}}\sin(\theta_{i_{k-1}})\cos(\theta_{i_{k-1}+1})\ldots\cos(\theta_{i_{k}-1})\sin(\theta_{i_k})\cos(\theta_{i_k+1})\cos(\theta_{i_k+2})\ldots\cos(\theta_r).$$
As each off-diagonal term switches the index, there must be an even number of off-digonal terms in the trace and therefore $|I|$ is even.
Then by induction
$$ (A(u_1),\ldots,A(u_r)\gamma)^1_1 =  (-1)^{s(I)}\sin(\theta_I)\cos(\theta_{\hat{I}})e^{L/2-L_I}$$
where 
$$s(I) = (i_2-i_1 +1) + (i_4-i_3+1) +\ldots (i_k-i_{k-1} + 1)$$
and 
 $$L_I= \sum_{j=1}^{k} (-1)^{j}l_{i_j} = -l_{i_1}+l_{i_2}-l_{i_3}-\ldots + (-1)^{k}l_{i_k}$$

Similarly 
$$\frac{(A(u_1),\ldots,A(u_r))^2_2(I)}{(A(u_{1})\ldots A(u_{i_{k-2}}))^2_2(I_{k-1})}= $$
$$= (-e^{-l_{i_{k-1}}}\sin(\theta_{i_{k-1}}))(\cos(\theta_{i_{k-1}+1}))(\cos(\theta_{i_{k-1}+2}))\ldots(\cos(\theta_{i_k-1}))\left(-e^{l_{i_k}}\sin(\theta_{i_k})\right)(-\cos(\theta_{i_k+1}))\ldots(-\cos(\theta_{r}))
$$
$$ =(-1)^{r-i_{k}+2}e^{-l_{i_{k-1}}+l_{i_k}}\sin{\theta_{i_{k-1}}\cos(\theta_{i_{k-1}+1})\ldots\cos(\theta_{i_{k}-1})\sin(\theta_{i_k})\cos(\theta_{i_k+1})\cos(\theta_{i_k+2})\ldots\cos(\theta_r}).$$

Counting negative signs we have  $r-s(I) + |I|$ negative signs. 
$$ (A(u_1),\ldots,A(u_r)\gamma)^2_2 =  (-1)^{r-s(I)+|I|}\sin(\theta_I)\cos(\theta_{\hat{I}})e^{L_I-L/2}$$
As $|I|$  is even we get
$$ (A(u_1),\ldots,A(u_r)\gamma)^2_2 =  (-1)^{r+s(I)}\sin(\theta_I)\cos(\theta_{\hat{I}})e^{L_I-L/2}$$
giving the result.
\eproof

\section{Some examples}

We have from the calculations in the last section that
$$F_0(L) = 2\cosh(L/2)\qquad F_1(u,L) = 2\sinh(L/2)\cos\theta $$
$$ F_2(u_1,u_2,L) = 2(\cos\t_{1} \cos\t_{2}\cosh(L/2) + \sin\t_{1}\sin\t_{2}\cosh(L/2-l_{12}))$$
Calculating $F_3$ we have
$$F_3(u_1,u_2,u_3,L) = 2\sinh(L/2)\cos(\t_1) \cos(\t_2) \cos(\t_3)  + 2\sinh(L/2-l_{12})\sin(\t_1) \sin(\t_2)\cos(\t_3) $$
$$
  - 2\sinh(L/2-l_{13})  \sin(\t_1)\cos(\t_2)\sin(\t_3) +2\sinh(L/2-l_{23})  \cos(\t_1)\sin(\t_2) \sin(\t_3) $$
  
  Therefore we have
$$G_0(L) = 2\cosh(L/2)\qquad \qquad G_1(u_1,\ldots,u_n,L) =  2\sinh(L/2)\sum_{i=1}^n \cos{\theta_i} $$
$$G_2(u_1,\ldots,u_n,L) = 2\sum_{\stackrel{i,j=1}{i < j}}^n(\cos\t_{i} \cos\t_{j}\cosh(L/2) + \sin\t_{i}\sin\t_{j}\cosh(L/2-l_{ij}))$$
$$G_3(u_1,\ldots,u_n,L) = 2\sum_{i<j <k}\left(\sinh(L/2)\cos(\t_i) \cos(\t_j) \cos(\t_k)  + \sinh(L/2-l_{ij})\sin(\t_i) \sin(\t_j)\cos(\t_k) \right.$$
$$\left.
  - \sinh(L/2-l_{ik})  \sin(\t_i)\cos(\t_j)\sin(\t_k) +\sinh(L/2-l_{jk})  \cos(\t_i)\sin(\t_j) \sin(\t_k)\right) $$

As the functions $G_r$ do not depend on $k$, once we've calculated all derivatives less than $k$, we only need calculate $G_k$ to find the $k$th derivative.

For $k=3$ we have
$$t_\beta^{3}T_\gamma = \frac{1}{8} \left( B_{n,3,1}G_{1}(u_1,\ldots,u_n,L_\gamma) + B_{n,3,3}G_{3}(u_1,\ldots,u_n,L_\gamma)\right)$$
$$B_{n,3,3} = 3! = 6\qquad \qquad B_{n,3,1} = (n-1)\left( \begin{array}{c}
3\\
1,2\end{array}
\right) + \left( \begin{array}{c}
3\\
3\end{array}
\right) = 3(n-1)+1 = 3n-2$$

$$t_\beta^{3}T_\gamma= \frac{1}{8} ((6n-4)\sinh(L_\gamma/2)\sum_{i=1}^n \cos(\theta_i) +12\left(\sum_{i<j <k}\sinh(L_\gamma/2)\cos(\t_i) \cos(\t_j) \cos(\t_k) \right.$$
$$  + \sinh(L_\gamma/2 -l_{ij})\sin(\t_i) \sin(\t_j)\cos(\t_k)  - \sinh(L_\gamma/2-l_{ik})  \sin(\t_i)\cos(\t_j)\sin(\t_k) $$
$$\left.+\sinh(L_\gamma/2-l_{jk})  \cos(\t_i)\sin(\t_j) \sin(\t_k)\right)$$

\end{document}